\documentclass[10pts]{article}
\usepackage{amssymb}
\usepackage{amsfonts}
\usepackage{amsmath}

\begin{document}
\title{Volterra-Choquet nonlinear operators}
\author{Sorin G. Gal \\
Department of Mathematics and Computer Science, \\
University of Oradea, \\
Universitatii Street No.1, 410087, Oradea, Romania\\
E-mail: \textit{galsorin23@gmail.com}\\
and \\
Academy of Romanian Scientists, \\
Splaiul Independentei nr. 54\\
050094 Bucharest, Romania}

\date{}
\maketitle
\begin{abstract}
In this paper we study to what extend some properties of the classical linear Volterra operators could be transferred to the nonlinear Volterra-Choquet operators, obtained by replacing the classical linear integral with respect to the Lebesgue measure, by the nonlinear Choquet integral with respect to a nonadditive set function. Compactness, Lipschitz and cyclicity properties are studied.
\end{abstract}

\section{Introduction}

Inspired by the electrostatic capacity, G. Choquet has introduced in \cite{Choquet}
(see also \cite{Ch1986}) a concept of integral with respect to a non-additive set
function which, in the case when the underlying set function is a $\sigma
$-additive measure, coincides with the Lebesgue integral.

Choquet integral is proved to be a powerful and useful tool in decision making under risk
and uncertainty, finance, economics, insurance, pattern recognition, etc (see, e.g.,
\cite{WK1} and \cite{WY} as well as the references therein).

Many new interesting results were obtained as analogs in the framework of
Choquet integral of certain known results for the Lebesgue integral. In this
sense, we can mention here, for example, the contributions to function spaces theory in \cite{Cerda1},
to potential theory in \cite{Adams},
to approximation theory in \cite{Gal1}%
-\cite{Gal6} and to integral equations theory in \cite{Gal9}, \cite{Gal10}.

Now, for $1\le p <+\infty$  denoting $L^{p}[0, 1]=\{f:[0, 1]\to \mathbb{R} ; (L)\int_{0}^{1}|f(x)|^{p}d x <+\infty\}$, where the $(L)$ integral is that with respect to the Lebesgue measure, it is well-known that the classical Volterra linear operator introduced in 1896, is defined usually on $L^{2}[0, 1]$ by
\begin{equation}\label{Volt}
O(f)(x)=(L)\int_{0}^{x}f(t) dt, x\in [0, 1].
\end{equation}
Volterra operator has been studied and continues to be studied by many authors. The norm of
Volterra operator is $2/\pi$ (see the book \cite{Ha}, Problem 149). The Halmos' book also contains several nice results
related with Volterra operator. The asymptotic behaviour of the norm $\|V^{n}\|$
is described in \cite{LR}.  An interesting fact about the Volterra operator
is the determination of its invariant subspace lattice (see \cite{Co}, Chapter 4 and \cite{Br}, \cite{Dix}, \cite{Don}, \cite{Ka} and \cite{Sar}).
Compactness and cyclicity properties were studied in, e.g.,  \cite{Sa}, \cite{Gallardo}, \cite{Leon}, \cite{Leon2}.
Very recent papers on various other aspects of the Volterrs operator are, e.g.,  \cite{Liang}, \cite{Khad}, \cite{Baksi}, \cite{ter}, \cite{Lefevre},
to mention only a few.
Note that there is also a huge literature dealing with the Volterra operator in complex setting, but this aspect is out of the discussions in the present paper.

Let ${\cal{C}}$ be a $\sigma$-algebra of subsets in $[0, 1]$ and $\mu:{\cal{C}}\to [0, +\infty]$ be a monotone set function, i.e. satisfying $\mu(\emptyset)=0$ and $\mu(A)\le \mu(B)$ for all $A, B\in {\cal{C}}$, with $A\subset B$.

The goal of the present paper is to study the possibilities of extension of the properties of classical Volterra linear operator, to the so-called Volterra-Choquet
operator obtained by replacing the classical linear integral by the nonlinear Choquet integral, that is defined by
\begin{equation}\label{C-Volt}
V(f)(x)=(C)\int_{0}^{x}f(t) d\mu(t),
\end{equation} where
$\mu$ is a set function not necessarily additive.

The plan of the paper goes as follows.
Section 2 contains preliminaries on the Choquet integral and Section 3 presents a few general preliminaries on compactness of nonlinear operators.
In Section 4 we prove some compactness properties while in Section 5 we obtain some Lipschitz properties, for the Volterra-Choquet operators.
Section 6 presents cyclicity properties for a Volterra-Choquet operator with respect to a particular distorted Lebesgue measure.

\section{Preliminaries on Choquet integral}

Some known concepts and results concerning the Choquet integral can be summarized by the following.

{\bf Definition 2.1.} Suppose $\Omega\not=\emptyset$ and let ${\cal{C}}$ be a $\sigma$-algebra of subsets in $\Omega$.

(i) (see, e.g., \cite{WK1}, p. 63) The set function  $\mu:{\cal{C}}\to [0, +\infty]$ is called a monotone set function (or capacity) if $\mu(\emptyset)=0$ and $\mu(A)\le \mu(B)$ for all $A, B\in {\cal{C}}$, with $A\subset B$. Also,
$\mu$ is called submodular if
$$\mu(A\bigcup B)+\mu(A\bigcap B)\le \mu(A)+\mu(B), \mbox{  for all } A, B\in {\cal{C}}.$$
$\mu$ is called bounded if $\mu(\Omega)<+\infty$ and normalized if $\mu(\Omega)=1$.

(ii) (see, e.g., \cite{WK1}, p. 233, or \cite{Choquet}) If $\mu$ is a monotone set function on ${\cal{C}}$
and if $f:\Omega \to \mathbb{R}$ is ${\cal{C}}$-measurable (that is, for any Borel subset $B\subset \mathbb{R}$ it follows $f^{-1}(B)\in {\cal{C}}$), then for any $A\in {\cal{C}}$, the concept of Choquet integral is defined by
$$(C)\int_{A} f d\mu=\int_{0}^{+\infty}\mu\left (F_{\beta}(f)\bigcap A\right )d\beta+\int_{-\infty}^{0}\left [\mu\left (F_{\beta}(f)\bigcap A\right )-\mu(A)\right ]d \beta,$$
where we used the notation $F_{\beta}(f)=\{\omega\in \Omega; f(\omega)\ge \beta\}$.
Notice that if $f\ge 0$ on $A$, then in the above formula we get $\int_{-\infty}^{0}=0$.

The function $f$ will be called Choquet integrable on $A$ if $(C)\int_{A}f d\mu\in \mathbb{R}$.

(iii) (see, e.g., \cite{WK1}, p. 40) We say that the set function $\mu:{\cal{C}}\to [0, +\infty]$ is continuous from below, if for any sequence $A_{k}\in {\cal{C}}$, $A_{k}\subset A_{k+1}$, for all $k=1, 2, ..., $ we have $\lim_{k\to \infty}\mu(A_{k})=\mu(A)$, where $A=\bigcup_{k=1}^{\infty}A_{k}$.

Also, we say that $\mu$ is continuous from above, if for any sequence $A_{k}\in {\cal{C}}$, $A_{k+1}\subset A_{k}$, for all $k=1, 2, ..., $ we have $\lim_{k\to \infty}\mu(A_{k})=\mu(A)$, where $A=\bigcap_{k=1}^{\infty}A_{k}$.

In what follows, we list some known properties of the Choquet integral.

{\bf Remark 2.2.} If $\mu:{\cal{C}}\to [0, +\infty]$ is a monotone set function, then the following properties hold :

(i) For all $a\ge 0$ we have $(C)\int_{A}af d\mu = a\cdot (C)\int_{A}f d\mu$ (if $f\ge 0$ then see, e.g., \cite{WK1}, Theorem 11.2, (5), p. 228 and if $f$ is of arbitrary sign, then see, e.g., \cite{Denn}, p. 64, Proposition 5.1, (ii)).

(ii) In general (that is if $\mu$ is only monotone), the Choquet integral is not linear, i.e. $(C)\int_{A}(f+g)d\mu \not=(C)\int_{A}fd\mu + (C)\int_{A}gd\mu$.

In particular, for all $c\in \mathbb{R}$ and $f$ of arbitrary sign, we have (see, e.g., \cite{WK1}, pp. 232-233, or \cite{Denn}, p. 65) $(C)\int_{A}(f+c)d \mu = (C)\int_{A}f d\mu + c\cdot \mu(A)$.

If $\mu$ is submodular too, then for all $f, g$ of arbitrary sign and lower bounded, the property of subadditivity holds
(see, e.g., \cite{Denn}, p. 75, Theorem 6.3)
$$(C)\int_{A}(f + g) d\mu \le (C)\int_{A}f d\mu + (C)\int_{A}g d\mu.$$

However, in particular, the comonotonic additivity holds, that is if $\mu$ is a monotone set function and $f, g$ are ${\cal{C}}$-measurable and comonotone on $A$ (that is $(f(\omega)-f(\omega^{\prime}))\cdot (g(\omega)-g(\omega^{\prime}))\ge 0$, for all $\omega, \omega^{\prime}\in A$), then by, e.g., Proposition 5.1, (vi), p. 65 in \cite{Denn}, we have
$$(C)\int_{A}(f + g) d\mu = (C)\int_{A}f d\mu + (C)\int_{A}g d\mu.$$

(iii) If $f\le g$ on $A$ then
$(C)\int_{A}f d\mu \le (C)\int_{A}g d\mu$ (see, e.g., \cite{WK1}, p. 228, Theorem 11.2, (3) if $f, g\ge 0$ and p. 232 if $f, g$ are of arbitrary sign).

(iv) Let $f\ge 0$. If $A\subset B$ then $(C)\int_{A}f d \mu \le (C)\int_{B}f d\mu.$
In addition, if $\mu$ is finitely subadditive (that is, $\mu(\bigcup_{k=1}^{n} A_{k})\le \sum_{k=1}^{n}\mu(A_{k})$, for all $n\in \mathbb{N}$), then
$$(C)\int_{A\cup B}f d\mu \le (C)\int_{A}f d\mu + (C)\int_{B}f d\mu.$$

(v) It is immediate that $(C)\int_{A}1\cdot d\mu(t)=\mu(A)$.

(vi) The formula $\mu(A)=\gamma(m(A))$, where
$\gamma :[0, m(\Omega)]\to \mathbb{R}$ is an increasing and concave function, with $\gamma(0)=0$ and
$m$ is a bounded measure (or bounded but only finitely additive) on a $\sigma$-algebra on $\Omega$ (that is, $m(\emptyset)=0$ and $m$ is countably additive), gives simple examples of monotone and submodular set functions (see, e.g., \cite{Denn}, pp. 16-17). Such of set functions $\mu$ are also called distorsions of countably additive measures (or distorted measures).

If $\Omega=[a, b]$, then for the Lebesgue (or any Borel) measure $m$ on $[a, b]$, $\mu(A)=\gamma(m(A))$ give simple examples of bounded, monotone and submodular set functions on $[a, b]$.

In addition, if we suppose that $\gamma$ is continuous at $0$ and at $m([a, b])$, then by the continuity of $\gamma$ on the whole interval $[0, m([a, b])]$ and from the continuity from below of any Borel measure, it easily follows that the corresponding distorted measure also is continuous from below.

For simple examples, we can take $\gamma(t)=t^{p}$ with $0<p<1$, $\gamma(t)=\frac{2 t}{1+t}$, $\gamma(t)=1-e^{-t}$, $\gamma(t)=\ln(1+t)$ for $t\ge 0$
and $\gamma(t)=\sin(t/2)$ for $t\in [0, \pi]$.

Now, let us consider that in  the above definition of a distorted Lebesgue measure, $\mu(A)=\gamma(m(A))$, in addition $\gamma$ is considered strictly increasing and differentiable. In this case, if $f$ is nonnegative, nondecreasing and continuous, then (see, e.g., \cite{Su1}, Theorem I)
$$(C)\int_{0}^{x}f(s)d \mu(s)=\int_{0}^{x}\gamma^{\prime}(x-s)f(s)d s,$$
while if $f$ is nonnegative, nonincreasing and continuous, then (see, e.g., \cite{Su1}, Theorem A.1)
$$(C)\int_{0}^{x}f(s)d \mu(s)=\int_{0}^{x}\gamma^{\prime}(s)f(s)d s.$$
(vii) If $\mu$ is a countably additive bounded measure, then the Choquet integral $(C)\int_{A}f d\mu$ reduces to the usual Lebesgue type integral (see, e.g., \cite{Denn}, p. 62, or \cite{WK1}, p. 226).

(viii) Let ${\cal{C}}$ be a $\sigma$-algebra of subsets in $[0, 1]$ and $\mu:{\cal{C}}\to [0, +\infty]$ be a monotone set function.The analogs of the Lebesgue spaces in the context of capacities can be
introduced for $1\leq p<+\infty$ via the formulas
\[
{\mathcal{L}}_{p, \mu}([0, 1])=\{f:[0, 1]\to \mathbb{R} ; f \mbox{ is } {\cal{C}} \mbox{- measurable } \text{ and
}(C)\int_{0}^{1}|f(t)|^{p}\mathrm{d}\mu(t)<+\infty\}.
\]
When $\mu$ is a subadditive capacity (in particular, when $\mu$ is submodular), the functionals $\Vert\cdot\Vert_{{\mathcal{L}}_{p, \mu}}$ given by
\begin{align*}
\Vert f\Vert_{{\mathcal{L}}_{p, \mu}} &  =\left(  (C)\int_{\Omega
}|f(t)|^{p}\mathrm{d}\mu(t)\right)  ^{1/p}
\end{align*}
satisfy the triangle inequality too (see, e.g. Theorem 2, p. 5 in \cite{Cerda1} or Proposition 9.4, p. 109 in \cite{Denn} or
Theorems 3.5 and 3.7 in \cite{RSWang} or the comments in the  proof of  Theorem 3.4, Step 3 in \cite{Gal}.

Denoting
$${\mathcal{N}_{p}}=\{f\in {\mathcal{L}}_{p, \mu}([0, 1]); \left ((C)\int_{0}^{1}|f(t)|^{p}\mathrm{d}\mu(t)\right)^{1/p}=0\},$$
if $\mu$ is a submodular capacity, then the functionals $\Vert\cdot\Vert_{L_{p, \mu}%
([0, }$ given by
\begin{align*}
\Vert f\Vert_{L_{p, \mu}([0, 1])} &  =\left(  (C)\int_{0}^{1}
|f(t)|^{p}\mathrm{d}\mu(t)\right)  ^{1/p}
\end{align*}
satisfy the axioms of a norm on the quotient space $L_{p, \mu}([0, 1])={\mathcal{L}}_{p, \mu}([0, 1])/{\mathcal{N}_{p}}$ (see \cite{Denn}, p. 109, Proposition 9.4 for $p=1$ and p. 115 for arbitrary $p\ge 1$). If, in addition, $\mu$ is continuous  from below, then
$L_{p, \mu}([0, 1])$ is a  Banach space (see \cite{Denn}, pp. 11-12, Proposition 9.5) and $h\in {\mathcal{N}_{p}}$ if and only if $h=0$, $\mu$-a.e., meaning that there exists $N$ with $\mu^{*}(N)=0$, such that $h(\omega)=0$, for all $\omega\in \Omega \setminus N$ (see \cite{Denn}, p. 107, Corollary 9.2 and pp. 107-108 ). Here $\mu^{*}$ is he outer measure attached to $\mu$, given by the formula $\mu^{*}(A)=\inf\{\mu(B); A\subset B, B\in {\cal{C}}\}$.

Also, denote
$$L_{p, \mu}^{+}([0, 1])=\{f:[0, 1]\to \mathbb{R} ; f \mbox{ is } {\cal{C}} \mbox{-measurable and } (C)\int_{0}^{1}|f(t)|^{p}d\mu(t)<+\infty\}.$$

\section{Preliminaries on nonlinear compact operators}

In this section we present a few well-known concepts and general results on compactness of nonlinear operators which we need for the next sections.

{\bf Definition 3.1.} Let $A : X \to Y$ be a nonlinear operator between two metric spaces. We recall that $A$ is said to be compact
if it is continuous on $X$ and for any bounded $M\subset X$,  $A(M)$ is relatively compact in $Y$ (that is the closure $\overline{A(M)}$ is compact in $Y$). Recall here that a set $M\subset X$ is called bounded in the metric space $(X, \rho)$ if $diam(M)=\sup\{\rho(x, y) ; x, y\in M\}<+\infty$.

{\bf Remark 3.2} If $X$ and $Y$ are two normed spaces over $\mathbb{R}$ and $A$ is positive homogeneous, then it is easy to see that $A$ is compact if and only if it is continuous and $A[B_1]$ is relatively compact in $Y$, where $B_{r}=\{x\in X ; \|x\|\le r\}$. Indeed, let $M\subset X$ be bounded, that is there exists $r>0$, such that $M\subset B_r$. It is immediate that $\overline{A[M]}\subset \overline{A[B_{r}]}=r \overline{A(B_{1}]}$. Since $r \overline{A[B_{1}]}$ is compact, it is clear that the closed set $\overline{A[M]}$ is compact.

We also recall that a fundamental result in the study of the algebra of continuous functions on a compact Hausdorff space is the well-known Arzel\`a-Ascoli theorem, which can be stated as follows.

{\bf Theorem 3.3.} (\cite{Dunford}, IV.6.7.) {\it Let $X$ be a compact Hausdorff space and denote
$$C(X; \mathbb{R})=\{f:X\to \mathbb{R} ; f \mbox{ is continuous on } X\}.$$
Then a subset $F$ of $C(X; \mathbb{R})$ is relatively compact in the topology induced by the uniform norm, if and only if it is equicontinuous and pointwise bounded. Here pointwise bounded means that for any $x\in X$ we have $\sup\{|f(x)| ; f\in F\}<\infty$.}

\section{Compactness of the Volterra-Choquet operators}

This section contains some important properties of compactness for the Volterra-Choquet operators.

In this sense, firstly we need the following.

{\bf Theorem 4.1.} {\it Let $\mu$ be a monotone, submodular and continuous from below set function on all Borelian subsets in $[0, 1]$ (or all the Lebesgue measurable subsets in $[0, 1]$), $1 < p<+\infty$ and $1/p + 1/q=1$. Then,
for all $x, y\in [0, 1]$ with $x \le y$, and $f\in L_{p, \mu}^{+}([0, 1])$,
the Volterra Choquet operator $V(f)(x)=(C)\int_{[0, x]}f(t)d\mu(t)$ has the property
$$|V(f)(x)-V(f)(y)|\le \|f\|_{L_{p, \mu}([0, 1])}\cdot \mu([x, y])^{1/q}.$$}
{\bf Proof.} Indeed, without loss of generality, we may suppose that $x < y$. The submodularity of $\mu$ evidently implies the finitely subadditivity. Since $[0, y]=[0, x]\bigcup [x, y]$, by Remark 2.2, (iv), we obtain $V(f)(y)\le V(f)(x)+(C)\int_{x}^{y}f(t)d \mu(t)$, which by applying the H\"older's inequality too ((see Remark 2.2, (viii)) e.g., implies
$$0\le V(f)(y)-V(f)(x)\le (C)\int_{[x, y]}f(t)d \mu(t)=(C)\int_{[x, y]}f(t) \cdot 1\cdot  d\mu(t)$$
$$\le \|f\|_{L_{p, \mu}([0, 1])}\cdot \left ( (C)\int_{[x, y]}1 \cdot d\mu(t)\right )^{1/q}=\|f\|_{L_{p, \mu}([0, 1])}\cdot \mu([x, y])^{1/q}.$$
If $y \le x$, then in the statement we obtain a similar inequality, by replacing $\mu([x, y])$ with $\mu([y, x])$,
which proves the theorem. $\hfill \square$

{\bf Corollary 4.2.} {\it Suppose that $\mu$ is a distorted Lebesgue measure, that is $\mu(A)=\gamma(m(A))$, where $m$ is the Lebesgue measure and $\gamma:[0, 1]\to \mathbb{R}$ is nondecreasing, concave, continuous on $[0, 1]$ and $\gamma(0)=0$. For $1< p <+\infty$, denote $B_{p, \mu, 1}^{+}=\{f\in L_{p, \mu}^{+}([0, 1]) : \|f\|_{L_{p, \mu}([0, 1])}\le 1\}$. Then, $V(B_{p, \mu, 1}^{+})$ is an equicontinuous and uniformly bounded set of continuous functions defined on $[0, 1]$.}

{\bf Proof.} Since by Remark 2.2, (vi), any distorted Lebesgue set function is submodular and continuous from below (in fact, from above too), by Theorem 4.1, we obtain
$$|V(f)(x)-V(f)(y)|\le \gamma(|x-y|)^{1/q}, \mbox{ for all } x,  y\in [0, 1] \mbox{ and }f\in B_{p, \mu, 1}^{+}([0, 1]).$$
Let $\varepsilon >0$ be fixed. By the continuity of $\gamma$, there exists a $\delta>0$ (depending of course only on $\varepsilon$ and $\gamma$ and indepepndent of $f$), such that $\gamma(|x-y|)^{1/q}<\varepsilon$, for all $|x-y| < \delta$. This immediately implies the equicontinuity of the set of continuous functions $V(B_{p, \mu, 1}^{+})$.

Also, choosing $y=0$ in the above inequality, we obtain
$$|V(f)(x)|\le \gamma(x)^{1/q}\le \gamma(1)^{1/q}, \mbox{ for all } x\in [0, 1], f\in B_{p, \mu, 1}^{+},$$
which proves that  $V(B_{p, \mu, 1}^{+})$ is uniformly bounded.
$\hfill \square$

By using Definition 3.1, we can state the following.

{\bf Corollary 4.3.} {\it Under the hypothesis of Corollary 4.2, the Volterra-Choquet operator $V:L_{p, \mu}^{+}([0, 1])\to C_{+}[0, 1]\subset C[0, 1]$ is a nonlinear compact  operator.

Here $L_{p, \mu}^{+}([0, 1])$ is endowed with the metric generated by the $L_{p, \mu}$-norm in Remark 2.2, (viii)
and $C_{+}[0, 1]$ denotes the space of all nonnegative real-valued continuous functions on $[0, 1]$,
which is a metric space endowed with the metric generated by the uniform norm.}

{\bf Proof.} By Arzel\'a-Ascoli result in Theorem 3.3 and by Corollary 4.2, it follows that $\overline{V(B_{p, \mu, 1}^{+})}$ is compact.

Let $M\subset L_{p, \mu}^{+}([0, 1])$ be bounded, that is $d=diam(M)<+\infty$. For a fixed $x_{0}\in M$ and an arbitrary $x\in M$, we get
$$|\quad \|x\|_{L_{p, \mu}([0, 1])}-\|x_{0}\|_{L_{p, \mu}([0, 1])}\quad |\le \|x-x_{0}\|_{L_{p, \mu}([0, 1])}\le d,$$ which immediately implies
$\|x\|_{L_{p, \mu}([0, 1])}\le \|x_{0}\|_{L_{p, \mu}([0, 1])}+d$, that is
$$M\subset B^{+}_{p, \mu, \|x_{0}\|_{L_{p, \mu}([0, 1])}+d}=(\|x_{0}\|_{L_{p, \mu}([0, 1])}+d) B_{p, \mu, 1}^{+}.$$
Applying  $V$ and taking into account that by Remark 2.2, (i), $V$ is positive homogeneous, we get $V(M)\subset (\|x_{0}\|_{L_{p, \mu}([0, 1])}+d) V(B_{p, \mu, 1}^{+})$,
which implies $\overline{V(M)}\subset (\|x_{0}\|_{L_{p, \mu}([0, 1])}+d) \overline{V(B_{p, \mu, 1}^{+})}$.

Now, since in a metric space any closed subset of a compact set also is compact, it implies that $\overline{V(M)}$ is compact.

It remains to prove the continuity of the operator $V:L_{p, \mu}^{+}([0, 1])\to C_{+}[0, 1]$. For this purpose, let $f, g\in L_{p, \mu}([0, 1])$ (in fact, not necessarily nonnegative). From H\"older's inequality it is immediate that $L_{p, \mu}^{+}([0, 1])\subset L_{1, \mu}^{+}([0, 1])$. Since according to Remark 2.2, (ii), the Choquet integral is subadditive, by
$f(s)\le g(s)+|f(s)-g(s)|$, for all $s\in [0, 1]$, it follows that
$$V(f)(t)\le (C)\int_{0}^{t}g(s)d \mu(s)+(C)\int_{0}^{1}|f(s)-g(s)|d \mu(s).$$
This implies $V(f)(t)\le V(g)(t)+V(|f-g|)(t)$, for all $t\in [0, 1]$.

Also, by $g(s)\le f(s)+|g(s)-f(s)|$, for all $s\in [0, 1]$, by similar reasoning we obtain
$V(g)(t)\le V(f)(t)+V(|f-g|)(t)$, for all $t\in [0, 1]$, which combined with the above inequality, leads to the inequality valid for all $t\in [0, 1]$
\begin{equation*}
|V(f)(t)-V(g)(t)|\le V(|f-g|)(t)\le (C)\int_{0}^{1}|f(s)-g(s)|\cdot 1 d\mu(s).
\end{equation*}
Passing to supremum after $t\in [0, 1]$ in the left hand-side and then, applying the H\"older's inequality to the right-hand side,
we easily arrive to
$$\|V(f)-V(g)\|_{C[0, 1]}\le \mu([0, 1])^{1/q}\cdot \|f-g\|_{L_{p, \mu}([0, 1])},$$
from which easily follows the continuity of $V$.

Concluding, by Definition 3.1 all the above mean the compactness of $V:L_{p, \mu}^{+}([0, 1])\to C_{+}[0, 1]$.
$\hfill \square$

{\bf Remark 4.4.} In the case of $p=1$, Corollary 4.3 does not hold in general. Indeed, it is known that even in the very particular case when $\gamma(t)=t$ (that is when $\mu$ one reduces to the classical Lebesgue measure), the equicontinuity fails.

\section{Lipschitz type properties and compactness}

In this section, firstly we prove Lipschitz properties of the nonlinear Volterra-Choquet operator $V$, on the whole spaces $C[0, 1]$ and $L_{p, \mu}([0, 1])$ with $1\le p<+\infty$.

{\bf Theorem 5.1.} {\it Let $\mu$ be a monotone, submodular and continuous from below and from above set function on the class of all Borelian (or alternatively, on the class of all Lebesgue measurable) subsets of $[0, 1]$.

(i) If $f\in L_{1, \mu}([0, 1])$ then $V(f)\in L_{1, \mu}([0, 1])$ and for all $f, g\in L_{1, \mu}([0, 1])$ we have
$$\|V(f)-V(g)\|_{L_{1, \mu}([0, 1])}\le \mu([0, 1])\cdot \|f-g\|_{L_{1, \mu}([0, 1])}.$$
(ii) If $f\in C[0, 1]$ then $V(f)\in C[0, 1]$ and for all $f, g\in C[0, 1]$ we have
$$\|V(f)-V(g)\|_{C[0, 1]}\le \mu([0, 1])\cdot \|f-g\|_{C[0, 1]},$$
where $\|\cdot \|_{C[0, 1]}$ denotes the uniform norm on $C[0, 1]$.

(iii) Let $1 < p < +\infty$. If $f\in L_{p, \mu}([0, 1])$ then $V(f)\in L_{p, \mu}([0, 1])$ and for all $f, g\in L_{p, \mu}([0, 1])$ we have
$$\|V(f)-V(g)\|_{L_{p, \mu}([0, 1])}\le \mu([0, 1])\cdot \|f-g\|_{L_{p, \mu}([0, 1])}.$$}
{\bf Proof.} (i) Firstly, we need to show that if $f\in L_{1, \mu}([0, 1])$, then $V(f)\in L_{1, \mu}([0, 1])$. In the case when $f\ge 0$, the proof is simple, because $V(f)(x)=(C)\int_{0}^{x}f(t)d\mu(t)$ is evidently a nondecreasing function on $[0, 1]$. Therefore, $V(f)(x)$ is a Borel (Lebesgue) measurable function.

Suppose now that $f$ is bounded and has negative values too, that is there exist $M^{\prime}<0$ and $M>0$, such that $M^{\prime}\le f(t)\le M$,  $t\in [0, 1]$. By Definition 2.1, (ii), we have
$$V(f)(x)$$
$$=\int_{0}^{+\infty}\mu(\{t\in [0, x] ; f(t)\ge \alpha\})d\alpha + \int_{-\infty}^{0}\left [\mu(\{t\in [0, x] ; f(t)\ge \alpha\})-\mu([0, x])\right ]d\alpha$$
$$=\int_{0}^{M}\mu(\{t\in [0, x] ; f(t)\ge \alpha\})d\alpha + \int_{M^{\prime}}^{0}\left [\mu(\{t\in [0, x] ; f(t)\ge \alpha\})-\mu([0, x])\right ]d\alpha$$
\begin{equation}\label{Bounded}
=\int_{0}^{M}\mu(\{t\in [0, x] ; f(t)\ge \alpha\})d\alpha + \int_{M^{\prime}}^{0}\mu(\{t\in [0, x] ; f(t)\ge \alpha\})d \alpha+
M^{\prime} \mu([0, x]).
\end{equation}
Therefore, $V(f)(x)$ is the sum of two nondecreasing functions with a nonincreasing one, all of them being Borel (Lebesgue) measurable, implying that
$V(f)(x)$ is Borel (Lebesgue) measurable too.

Suppose now that $f$ is unbounded and has negative values too. By the above formula, we can write $V(f)(x)=F(x)+G(x)$,
with
$$F(x)=\int_{0}^{+\infty}\mu(\{t\in [0, x] ; f(t)\ge \alpha\})d\alpha,$$
$$G(x)=\int_{-\infty}^{0}[\mu(\{t\in [0, x] ; f(t)\ge \alpha\})-\mu([0, x])]d\alpha$$
$$=\lim_{n\to +\infty}\int_{-n}^{0}[\mu(\{t\in [0, x] ; f(t)\ge \alpha\})-\mu([0, x])]d\alpha.$$
Evidently $F(x)$ is nondecreasing and therefore Borel (Lebesgue) measurable. Then, since for each $n\in \mathbb{N}$,
$$\int_{-n}^{0}[\mu(\{t\in [0, x] ; f(t)\ge \alpha\})-\mu([0, x])]d\alpha$$
$$=\int_{-n}^{0}\mu(\{t\in [0, x]; f(t)\ge \alpha\})d\alpha-n\mu([0, x])$$
is Borel (Lebesgue) measurable as a difference of two measurable functions, it follows that $G(x)$ is Borel (Lebesgue) measurable as a pointwise limit of Borel (Lebesgue) measurable functions.

Concluding, in this case too we have that $V(f)(x)$ is Borel (Lebesgue) measurable.

Then, we have
$$(C)\int_{0}^{1}|V(f)(x)|d\mu(x)=(C)\int_{0}^{1}\left |(C)\int_{0}^{x}f(t)d\mu(t)\right |d\mu(x)$$
$$\le (C)\int_{0}^{1}\left [(C)\int_{0}^{1}|f(t)|d\mu(t)\right ]d\mu(x)
=\|f\|_{L_{1, \mu}([0, 1])}\cdot (C)\int_{0}^{1}d\mu(x)$$
$$=\mu([0, 1])\cdot \|f\|_{L_{1, \mu}([0, 1])}<+\infty.$$
Let $f, g\in L_{1, \mu}([0, 1])$. Since according to Remark 2.2, (ii), the Choquet integral is subadditive, by
$f(s)\le g(s)+|f(s)-g(s)|$, for all $s\in [0, 1]$, it follows that
$$V(f)(t)=(C)\int_{0}^{t}f(s)d \mu(s)\le (C)\int_{0}^{t}g(s)d \mu(s)+(C)\int_{0}^{1}|f(s)-g(s)|d \mu(s),$$
which implies $V(f)(t)\le V(g)(t)+V(|f-g|)(t)$, for all $t\in [0, 1]$.

Also, by $g(s)\le f(s)+|g(s)-f(s)|$, for all $s\in [0, 1]$, by similar reasoning we obtain
$V(g)(t)\le V(f)(t)+V(|f-g|)(t)$, for all $t\in [0, 1]$. Combining these two inequalities, it is immediate that for all $t\in [0, 1]$ we have
\begin{equation}\label{Lipschitz}
|V(f)(t)-V(g)(t)|\le V(|f-g|)(t)\le (C)\int_{0}^{1}|f(s)-g(s)|d\mu(s)=\|f-g\|_{L_{1, \mu}([0, 1])}.
\end{equation}
Applying the Choquet integral, we obtain
$$\|V(f)-V(g)\|_{L_{1, \mu}([0, 1])}\le \mu([0, 1])\cdot \|f-g\|_{L_{1, \mu}([0, 1])}.$$
(ii)  Firstly, we need to show that if $f\in C[0, 1]$ then $V(f)\in C[0, 1]$. In the case when $f\ge 0$, the proof is as follows.

Let $x_{n}, x\in [0, 1]$ with $x_{n}\to x$.
Denoting $A_{n}(x, \alpha)=\{s\in [0, x_{n}]; f(s)\ge \alpha\}$, $A(x, \alpha)=\{s\in [0, x]; f(s)\ge \alpha\}$ and $M=\sup\{|f(x)|; x\in [0, 1]\}$, by the Definition 2.1, (ii) of the Choquet integral (for positive functions), we easily get
$$|V(f)(x_{n})-V(f)(x)|$$
$$=\left |\int_{0}^{M}\mu(\{s\in [0, x_{n}]; f(s)\ge \alpha\})d \alpha
-\int_{0}^{M}\mu(\{s\in [0, x]; f(s)\ge \alpha\})d \alpha \right |$$
$$=\left |\int_{0}^{M}\mu(A_{n}(x, \alpha))d \alpha - \int_{0}^{M}\mu(A(x, \alpha))d \alpha\right |.$$
According to, e.g., \cite{Siret}, p. 211, Exercise 5.5.21, it suffices to consider that $x_{n}\nearrow x$ and that $x_{n}\searrow x$.

If $x_{n}\nearrow x$, then since $A_{n}(x, \alpha)\subset A_{n+1}(x, \alpha)\subset A(x, \alpha)\subset [0, 1]$, for all $n\in \mathbb{N}$ and $\alpha \in [0, M]$, it follows $\mu(A_{n}(x, \alpha))\le \mu([0, 1])$, for all $n\in \mathbb{N}$, $\alpha\in [0, M]$ and by the continuity from below of $\mu$ we get
$\lim_{n\to \infty}\mu(A_{n}(x, \alpha))=\mu(A(x, \alpha))$, for all $\alpha\in [0, M]$. Passing then to limit under the integrals (which can be considered of Lebesgue type), we immediately obtain $\lim_{n\to \infty}|V(f)(x_{n})-V(f)(x)|=0$.

If $x_{n}\searrow x$, then since $A(x, \alpha)\subset A_{n+1}(x, \alpha)\subset A_{n}(x, \alpha)\subset [0, 1]$, for all $n\in \mathbb{N}$, $\alpha\in [0, M]$, it follows $\mu(A_{n}(x, \alpha))\le \mu([0, 1])$, for all $n\in \mathbb{N}$, $\alpha\in [0, M]$ and by the continuity from above of $\mu$ we get
$\lim_{n\to \infty}\mu(A_{n}(x, \alpha))=\mu(A(x, \alpha))$, for all $\alpha\in [0, M]$. Again, passing  to limit under the integrals, we immediately obtain $\lim_{n\to \infty}|V(f)(x_{n})-V(f)(x)|=0$.

In conclusion, in this case $V(f)\in C[0, 1]$.

Now, if $f\in C[0, 1]$ has negative values too, there exist $M^{\prime}<0$ and $M>0$, such that $M^{\prime}\le f(t)\le M$, for all $t\in [0, 1]$.
Then, by Definition 2.1, (ii) and by (\ref{Bounded}), we get
$$V(f)(x)=$$
$$\int_{0}^{M}\mu(\{t\in [0, x] ; f(t)\ge \alpha\})d\alpha + \int_{M^{\prime}}^{0}\mu(\{t\in [0, x] ; f(t)\ge \alpha\})d \alpha+
M^{\prime} \cdot \mu([0, x]).$$
Now, taking $x_{n}\to x$ and reasoning for each of the three terms on the right-hand side above as we reasoned in the case when $f\ge 0$, we easily obtain that $V(f)\in C[0, 1]$ in this case too.

The Lipschitz inequality, follows immediately by passing to supremum after $t\in [0, 1]$ in formula (\ref{Lipschitz}).

(iii) Let $f\in L_{p, \mu}([0, 1])$ with $1<p<+\infty$. Since from H\"older's inequality we get $L_{p, \mu}([0, 1])\subset L_{1, \mu}([0, 1])$, reasoning as at the point (i), it follows that $V(f)(x)$ is Borel (Lebesgue) measurable.

By (\ref{Lipschitz}) and using the H\"older's inequality (see Remark 2.2, (viii)), it follows
$$|V(f)(t)-V(g)(t)|\le \int_{0}^{1}|f(s)-g(s)|\cdot 1 \cdot d\mu(s)$$
$$\le \left ((C)\int_{0}^{1}|f(s)-g(s)|^{p}d\mu(s)\right )^{1/p}\cdot \mu([0, 1])^{1/q}.$$
Taking to the power $p$ both members of the above inequality, applying the Choquet ntegral on $[0, 1]$ with respect to $t$ and then taking to the power $1/p$,
we obtain
$$\|V(f)-V(g)\|_{L_{p, \mu}([0, 1])}\le \mu([0, 1])^{1/p}\cdot \mu([0, 1])^{1/q}\cdot \|f-g\|_{L_{p, \mu}([0, 1])}$$
$$=\mu([0, 1])\cdot \|f-g\|_{L_{p, \mu}([0, 1])},$$
which ends the proof of the theorem. $\hfill \square$

{\bf Remark 5.2.} Any distorted Lebesgue measure (defined as in Remark 2.2, (vi)) satisfies the hypothesis in Theorem 5.1.

{\bf Remark 5.3.} By Theorem 5.1, (i) and (iii), for the norm of  $V$ given by
$$\||V\||_{L_{p, \mu}([0, 1])}=\sup\left \{\frac{\|V(f)-V(g)\|_{L_{p, \mu}([0, 1])}}{\|f-g\|_{L_{p, \mu}([0, 1])}}; f, g\in L_{p, \mu}([0, 1]), f\not=g \right \}$$
we have $\||V\||_{L_{p, \mu}([0, 1])}\le \mu([0, 1])$, $1\le p <+\infty$.

As an application of Theorem 5.1, (ii), the following compactness property holds.

{\bf Corollary 5.4.} {\it  Suppose that $\mu$ is a distorted Lebesgue measure, that is $\mu(A)=\gamma(m(A))$, where $m$ is the Lebesgue measure, $\gamma:[0, 1]\to \mathbb{R}$ is nondecreasing, concave, continuous on $[0, 1]$ and $\gamma(0)=0$. Then, the Volterra-Choquet operator $V:C[0, 1]\to C[0, 1]$ defined by $V(f)(x)=(C)\int_{0}^{x}f(t)d \mu(t)$ is a nonlinear compact  operator.}

{\bf Proof.}  By Theorem 5.1, (ii), we have $V:C[0, 1]\to C[0, 1]$ and obviously that the Lipschitz property implies the continuity of $V$.

Let $B_{\mu, 1}=\{f\in C[0, 1] ; \|f\|_{C[0, 1]}\le 1\}$. If we prove that $V(B_{\mu, 1})$ is equicontinuous and pointwise bounded in $C[0, 1]$, then since by Remark 2.2, (i) $V$ is positive homogeneous, applying Remark 3.2
and Theorem 3.3, we will  immediately obtain that $V$ is compact.

Indeed, firstly for $f\in B_{\mu, 1}$ we get
$$|V(f)(x)|$$
$$=|(C)\int_{[0, x]}f(t)|d\mu(t)|\le \|f\|_{C[0, 1]}\cdot \gamma(m([0, 1]))=\gamma(1), \mbox{ for all } x\in [0, 1]$$
which implies $\|V(f)\|_{C[0, 1]}\le \gamma(1)$, for all $f\in B_{\mu, 1}$. This means that $V(B_{\mu, 1})$ is uniformly bounded (more than pointwise bounded).

Now, applying (\ref{Bounded}) for $M^{\prime}=-1$, $M=+1$, for all $f\in B_{\mu, 1}$ and $x\in [0, 1]$, we get
$$V(f)(x)=$$
$$=\int_{0}^{1}\mu(\{t\in [0, x] ; f(t)\ge \alpha\})d\alpha + \int_{-1}^{0}\mu(\{t\in [0, x] ; f(t)\ge \alpha\})d\alpha-\mu([0, x]).$$
Let $x_{n}\to x$. We easily get
$$|V(f)(x_{n})-V(f)(x)|$$
$$\le \int_{0}^{1}|\mu(\{t\in [0, x_{n}] ; f(t)\ge \alpha\})-\mu(\{t\in [0, x] ; f(t)\ge \alpha\})|d\alpha$$
$$+\int_{-1}^{0}|\mu(\{t\in [0, x_{n}] ; f(t)\ge \alpha\})-\mu(\{t\in [0, x] ; f(t)\ge \alpha\})|d\alpha$$
$$+|\mu([0, x_{n}])-\mu([0, x])|.$$
Firstly, suppose that $x_{n}\searrow x$. Since $[0, x_{n}]=[0, x]\cup [x, x_{n}]$ and $\mu$ submodular implies that $\mu$ is subadditive,
we obtain
$$0\le \mu(\{t\in [0, x_{n}] ; f(t)\ge \alpha\})-\mu(\{t\in [0, x] ; f(t)\ge \alpha\})$$
$$\le \mu(\{t\in [x, x_{n}] ; f(t)\ge \alpha\})\le \mu([x, x_{n}]),$$
which immediately implies
$$|V(f)(x_{n})-V(f)(x)|\le 3 \mu([x, x_{n}]).$$
The continuity of $\mu$ from above, immediately implies that $\lim_{n\to \infty}|V(f)(x_{n})-V(f)(x)|=0$, independent of $f$.

If $x_{n}\nearrow x$ we write $[0, x]=[0, x_{n}]\cup [x_{n}, x]$ and by analogous reasonings (since $\mu$ is continuous from below too) we get
$$|V(f)(x_{n})-V(f)(x)|\le 3 \mu([x_{n}, x]).$$
Concluding,  $V(B_{\mu, 1})$ is equicontinuous and therefore Corollary 5.4 is proved. $\hfill \square$

\section{Cyclicity}

Firstly recall the following known concepts.

{\bf Definition 6.1.} Let $(X, \|\cdot \|)$ be a Banach space on $K = \mathbb{R}$ (the real line) or $\mathbb{C}$
(the complex plane).

(i) The (not necessarily linear) continuous operator $T:X\to X$ is called cyclic, if there exists $x\in X$ such that the linear span of $Orb(T, x)$ is dense in $X$. Here $Orb(T, x)=\{x, T(x), T^{2}(x), ...,T^{n}(x), ...,  \}$.

(ii) $T$ is called hypercyclic, if there exists $x\in X$, such that the orbit
$Orb(T, x)$ is dense in $X$. Of course, if $X$ supports such an operator, then $X$ must be
separable.

(iii) $T$ is called supercyclic if there exists $x\in X$ such that the set ${\cal{M}}(x)=\{\lambda y ; y\in Orb(T, x); \lambda \in K\}$ is dense in $X$.

In the classical case, it is well-known that the Volterra operator $O$ given by (\ref{Volt}) and the identity plus Volterra operator, $I + O$, are cyclic operators on $L^{2}[0, 1]$, but they cannot be supercyclic and hypercyclic, see, e.g., \cite{Leon} and \cite{Gallardo}.

In what follows, we deal with cyclic type properties of the Volterra-Choquet operators with respect to a particular distorted Lebesgue measure. The problem of cyclic properties in the most general case remains open.

{\bf Theorem 6.2.} {\it Suppose that $\mu$ is the distorted Lebesgue measure given by $\mu(A)=\gamma(m(A))$, where $m$ is the Lebesgue measure, $\gamma:[0, 1]\to \mathbb{R}$ is $\gamma(x)=1-e^{-x}$ and that $f_{0}(x)=1$, for all $x\in [0, 1]$. Then, for the Volterra-Choquet operator $V$ with respect to $\mu$,  we have
$$Orb(V, f_{0})=\left \{1, \left \{1-e^{-x}\sum_{k=0}^{n-1}\frac{x^{k}}{k !}, n=1, 2, ...,\right \}\right \}.$$}

{\bf Proof.}  Firstly, by direct calculation we get
$$V(f_{0})(x)=\int_{0}^{\infty}\mu(\{s\in [0, x]; 1\ge \alpha\})d \alpha=\int_{0}^{1}\mu(\{s\in [0, x]; 1\ge \alpha\})d \alpha$$
$$=\mu([0, x])=\gamma(m([0, x]))=\gamma(x)=1-e^{-x}.$$
But, according to Proposition  1 in \cite{Su1}, see also Remark 2.2, (vi), (since $\gamma(x)$ is strictly increasing) and nonnegative), it follows
$$V^{2}(f_{0})(x)$$
$$=(C)\int_{0}^{x}V(f_{0})(s)d\mu(s)=\int_{0}^{x}\gamma^{\prime}(x-s)\gamma(s)d s$$
$$=\int_{0}^{x}e^{-(x-s)}(1-e^{-s})d s=1-e^{-x}-xe^{-x}=V(f_{0})(x)-xe^{-x}:=g_{1}(x).$$
Since $g_{1}$ is strictly increasing and nonnegative on $[0, 1]$ ($g_{1}(0)=0$), again by Proposition 1 in \cite{Su1}, it follows
$$V^{3}(f_{0})(x)=\int_{0}^{x}\gamma^{\prime}(x-s)V^{2}(f_{0})(s)d s=\int_{0}^{x}e^{-(x-s)}[1-e^{-s}-se^{-s}]ds$$
$$=e^{-x}\int_{0}^{x}(e^{s}-1-s)d s=e^{-x}(e^{x}-1-x-\frac{x^{2}}{2})=1-e^{-x}-xe^{-x}-e^{-x}\cdot \frac{x^{2}}{2}.$$
$$=V^{2}(f_{0})(x)-e^{-x}\cdot \frac{x^{2}}{2}:=g_{2}(x).$$
Since $g_{2}$ is strictly increasing and nonnegative on $[0, 1]$, we get
$$V^{4}(f_{0})(x)=\int_{0}^{x}\gamma^{\prime}(x-s)V^{3}(f_{0})(s)d s=\int_{0}^{x}e^{-(x-s)}[V^{2}(f_{0})(s)-e^{-s}\cdot \frac{s^{2}}{2}]d s$$
$$=\int_{0}^{x}e^{-(x-s)}V^{2}(f_{0})(s)d s-e^{-x}\int_{0}^{x}\frac{s^{2}}{2}d s=V^{2}(f_{0})(x)-e^{-x}\cdot \frac{x^{2}}{2}-e^{-x}\cdot \frac{x^{3}}{3 !}$$
$$=V^{3}(f_{0})(x)-e^{-x}\cdot \frac{x^{3}}{3 !}.$$
Continuing this kind of reasoning, we easily arrive at the general recurrence formula
$$V^{n}(f_{0})(x)=V^{n-1}(f_{0})(x)-e^{-x}\cdot \frac{x^{n-1}}{(n-1) !},$$
for all $n\ge 2$, which proves the theorem.
$\hfill \square$

{\bf Corollary 6.3.} {\it Let $V$ be the Volterra-Choquet operator with respect to the distorted Lebesgue measure $\mu$ in Theorem 6.2.

(i) As mapping $V:C[0, 1]\to C[0, 1]$, $V$ is a cyclic operator, with respect to the density induced by the uniform norm ;

(ii) Also, as mapping $V:L_{p, \mu}([0, 1])\to L_{p, \mu}([0, 1])$, $1\le p <+\infty$, $V$ is a cyclic operator with respect to the density induced by the $\|\cdot \|_{L_{p, \mu}([0, 1])}$-norm.}

{\bf Proof.} (i) By Theorem 6.2, it is immediate that span $Orb(V, f_{0})$ contains the countable subset given by  $\left \{1, e^{-x}\cdot \frac{x^{n}}{n !}, n=0, 1, 2, ....\right \}$, whose linear span is evidently dense in $C[0, 1]$ due to the Weierstrass approximation theorem by uniformly convergent sequences of polynomials.

(ii) It suffices to show that the set of polynomials is dense in $L_{p, \mu}([0, 1])$ with respect to the norm $\|\cdot \|_{L_{p, \mu}([0, 1])}$.
Indeed, since $\gamma^{\prime}(0)<+\infty$, by Remark 3.3. and Corollary 3.4 in \cite{Gal}, for the Bernstein-Durrmeyer-Choquet polynomials denoted by $D_{n, \mu}(f)(x)$ one has
$$\|D_{n, \mu}(f)-f\|_{L_{p, \mu}([0, 1])}\le c_{0}K(f; 1/\sqrt{n})_{L_{p, \mu}},$$
where $C^{1}_{+}([0, 1])$ denotes the the space of all differentiable, nonnegative functions with $g^{\prime}$ bounded on $[0, 1]$ and
$K(f; t)_{L_{p, \mu}}=\inf_{g\in C^{1}_{+}([0, 1])}\{\|f-g\|_{L_{p, \mu}([0, 1])}+t \|g^{\prime}\|_{C[0, 1]}\}$. According to Remark 3.3 in \cite{Gal}, it will suffices to prove that
$\lim_{t\to 0}K(f; t)_{L_{p, \mu}}=0$, for all $f\ge 0$.  Indeed, by e.g., \cite{Ivan}, there exist a sequence of non negative polynomials $(P_{n})_{n\in \mathbb{N}}$, such that $\|f-P_{n}\|_{p}\to 0$ as $n\to \infty$. For arbitrary $\varepsilon >0$, let $P_{m}$ be such
$\|f-P_{m}\|_{p}< \varepsilon/2$. Then for all $t \in (0,\, \varepsilon/(2 \|P^{\prime}_{m}\|_{\infty}))$, we get
$$
K(f; t)_{p}\ \le\ \|f-P_{m}\|_{p}\ +\ t \|P_{m}^{\prime}\|_{\infty}\ <\ \varepsilon/2\ +\ \varepsilon/2\ =\ \varepsilon,
$$
which proves our assertion. $\hfill \square$

{\bf Remark 6.4.} Using similar calculations with those in the proof of Theorem 6.2, we easily obtain that the operator $U=I+V$ satisfies the cyclicity properties in Corollary 6.3.

{\bf Remark 6.5.} The question that the Volterra-Choquet operator $V$ in Theorem 6.2 is, or is not, hypercyclic or supercyclic
remains unsettled. We observe that for $f_{0}=1$, neither $Orb(V, f_{0})$ and nor ${\cal{M}}(f_{0})$ are not dense in $C[0, 1]$.

{\bf Declaration of interest :} None.


\begin{thebibliography}{99}

\bibitem{Adams} D. R. Adams, Choquet integrals in potential theory,
Publicacions Matematiques, {\bf 42} (1998), 3-66.

\bibitem{Baksi} O. Baksi, T. Khan, J. Lang and V. Musil, Strict S-numbers of the Volterra operator, Proc. Amer. Math. Soc., {\bf 146} (2018), no. 2, 723-731.

\bibitem{Br} M. S. Brodskii, On a problem of I.M. Gelfand, Uspekhi Mat. Nak., {\bf 12} (1957), 129-
132.

\bibitem{Cerda1} J. Cerd\`{a}, J., Mart\'{\i}n and P. Silvestre,
Capacitary function spaces, Collect. Math., {\bf 62} (2011), 95-118.

\bibitem{Choquet} G. Choquet, Theory of capacities, Annales de l' Institut
Fourier, {\bf 5} (1954), 131-295.

\bibitem{Ch1986} G. Choquet, La naissance de la th\'{e}orie des capacit\'{e}s:
r\'{e}flexion sur une exp\'{e}rience personnelle, Comptes rendus de
l'Acad\'{e}mie des sciences. S\'{e}rie g\'{e}n\'{e}rale, La Vie des sciences,
{\bf 3} (1986), 385--397.

\bibitem{Co} J. B. Conway, A Course in Operator Theory, Rhode Island, Amer. Math. Soc., 2000.

\bibitem{Denn}
D. Denneberg, Non-Additive Measure and Integral, Kluwer Academic Publisher, Dordrecht, 1994.

\bibitem{Dix} J. Dixmier, Les operateurs permutables \`a l'operateur integral, Portugal Math., {\bf 8}
(1949), 73-84.

\bibitem{Don} W. F. Donoghue, The latice of invariant subspaces of quasi-nilpotent completely
continuous transformation, Pacific J. Math., {\bf 7} (1957), 1031-1035.

\bibitem{Dunford} N. Dunford and J. T. Schwartz, Linear Operators, vol. {\bf 1}, Wiley-Interscience, 1958.


\bibitem{ter} A. F. M. ter Elst and J. Zem\'anek, Contractive polynomials of the Volterra operator, Studia Math., {\bf 240} (2018), no. 3, 201-211.

\bibitem{Gal1} S. G. Gal and B. D. Opris, Uniform and pointwise
convergence of Bernstein-Durrmeyer operators with respect to monotone and
submodular set functions, J. Math. Anal. Appl., {\bf 424} (2015), 1374-1379.

\bibitem{Gal2} S. G. Gal, Approximation by Choquet integral
operators, Ann. Mat. Pura Appl., {\bf 195} (2016), no. 3, 881-896.

\bibitem{Gal} S. G. Gal and S. Trifa, Quantitative estimates in
$L^{p}$-approximation by Bernstein-Durr\-me\-yer-Cho\-qu\-et
operators with respect to distorted Borel measures, Results Math., {\bf 72} (2017),
no. 3, 1405-1415.

\bibitem{Gal4} S. G. Gal, Uniform and pointwise quantitative
approximation by Kantorovich-Choquet type integral operators with respect to
monotone and submodular set functions, Mediterr. J. Math., {\bf 14} (2017), no. 5,
art. 205, 12 pp.

\bibitem{Gal6} S. G. Gal, Quantitative approximation by nonlinear
Picard-Choquet, Gauss-Weierstrass-Choquet and Poisson-Cauchy-Choquet singular
integrals, Results Math., {\bf 73} (2018), no. 3, art. 92, 23 pp.

\bibitem{Gal9} S. G. Gal, Fredholm-Choquet integral equations, J.
Integral Equations and Applications,
https://projecteuclid.org/euclid.jiea/1542358961 (under press).

\bibitem{Gal10} S. G. Gal, Volterra-Choquet integral equations, J.
Integral Equations and Applications,
https://projecteuclid.org/euclid.jiea/1541668067 (under press).

\bibitem{Gallardo} E. Gallardo and A. Montes, Volterra operator is not supercyclic, Integral Equations
and Operator Theory, {\bf 50} (2004), no. 2, 211-216.

\bibitem{Ha} P. R. Halmos, A Hilbert Space Problem Book, New York, Springer Verlag, 1967.

\bibitem{Ivan} K. G. Ivanov, On a new characteristic  of functions, II. Direct and converse theorems
for the best algebraic approximation in $C[-1, 1]$ and $L^{p}[-1, 1]$, PLISKA Stud. Math. Bulgar., {\bf 5}, (1983), 151-163.

\bibitem{Ka} G. K. Kalisch, On similarity, reducing manifolds, and unitary equivalence of certain
Volterra operators, Ann. Math., {\bf 66} (1957), 481-494.

\bibitem{Khad} L. Khadkhuu and D. Tsedenbayar, A note about Volterra operator. Math. Slovaca, {\bf 68} (2018), no. 5, 1117-1120.

\bibitem{Khad2} L. Khadkhuu and  D. Tsedenbayar, On the numerical range and
numerical radius of the Volterra operator, The Bulletin of Irkutsk
State University. Series Mathematics, {\bf 24} (2018), 102-108.

\bibitem{Lefevre} P. Lef\`evre, The Volterra operator is finitely strictly singular from $L^{1}$ to $L_{\infty}$,
J. Approx. Theory {\bf 214} (2017), 1-8.

\bibitem{Leon} F. Le\'on-Saavedra and A. Piqueras-Lerena, Cyclic properties of Volterra operator,
Pacific J. Math., {\bf 211} (2003), no. 1, 157-162.

\bibitem{Leon2} F. Le\'on-Saavedra and A. Piqueras-Lerena, Super convex-cyclicity and the Volterra operator,
Monatsh. Math., {\bf 177} (2015), no. 2, 301-305.

\bibitem{Liang} Yu-Xia, Liang and Rongwei, Yang,  Energy functional of the Volterra operator,
Banach J. Math. Anal., {\bf 13} (2019), no. 2, 255-274.

\bibitem{LR} G. Little and J. B. Reade, Estimates for the norm of the n-th indefinite integral,
Bull. London Math. Soc., {\bf 30} (1998), 539-542.

\bibitem{Sa} N. H. Salas, Supercyclicity and weighted shifts, Studia Math., {\bf 135} (1999), 55-74.

\bibitem{Sar} D. Sarason, A remark on the Volterra operator, J. Math. Anal. Appl., {\bf 12} (1965),
244-246.

\bibitem{Siret} Gh. Siretchi, Differential and Integral Calculus, vol. II, Edit. Stiint. Enciclop., Bucharest, 1985.

\bibitem{Su1} M. Sugeno, A note on derivatives of functions with respect to fuzzy
measures, Fuzzy Sets Syst., {\bf 222} (2013), 1-17.

\bibitem{Ts} D. Tsedenbayar, On the power boundedness of certain Volterra operator pencils,
Studia Math., {\bf 156} (2003), 59-66.

\bibitem{RSWang}
R. S. Wang, Some inequalities and convergence theorems for Choquet integrals, J. Appl. Math. Comput.,
{\bf 35}(2011), 305-321.

\bibitem{WK1}
Z. Wang and G. J. Klir, Generalized Measure Theory, Springer, New York, 2009.

\bibitem {WY}Z. Wang and J.-A. Yan, A selective overview of applications of
Choquet integrals, Advanced Lectures in Mathematics, pp. 484--515, Springer, 2007.

\end{thebibliography}
\end{document}